\documentclass[12 pt, a4 paper]{article}
\usepackage{amssymb,amsmath,graphicx,enumerate}
\usepackage{longtable,mathrsfs}
\usepackage{hyperref}
\numberwithin{equation}{section}
\newtheorem{theorem}{Theorem}[section]

\newtheorem{definition}{Definition}[section]

\allowdisplaybreaks
\begin{document}
\begin{center}\begin{large}
New approach to Minkowski fractional inequalities using generalized k-fractional integral operator
\end{large}\end{center}
\begin{center}
                 $Vaijanath  \, L. Chinchane $

Department of Mathematics,\\
Deogiri Institute of Engineering and Management\\
Studies Aurangabad-431005, INDIA\\
chinchane85@gmail.com 
\end{center}
\begin{abstract}
 In this paper, we obtain new results related to Minkowski fractional integral inequality using generalized k-fractional integral operator which is in terms of the Gauss hypergeometric function.
\end{abstract}
\section{Introduction }
 \paragraph{} In the last decades many researchers have worked on fractional integral inequalities using Riemann-Liouville, generalized Riemann-Liouville, Hadamard and Siago, see \cite{C1,C2,C3,D1,D2,D3,D4}. W. Yang \cite{YA} proved the Chebyshev and Gr$\ddot{u}$ss-type integral inequalities for Saigo fractional integral operator. S. Mubeen and S. Iqbal \cite{MU} has proved the Gr$\ddot{u}$ss-type integral inequalities generalized k-fractional integral. In \cite{BA1,C5,KI2,YI} authors have studied some fractional integral inequalities using generalized k-fractional integral operator (in terms of the Gauss hypergeometric function). Recently many researchers have shown development of fractional integral inequalities associated with hypergeometric functions, see \cite{SH1,KI2,P1,R1,S1,SA,V1,W1,YI}.  Also, in \cite{C2,D1} authors established reverse Minkowski fractional integral inequality using Hadamard and Riemann-Liouville integral operator respectively.
\paragraph{}In literature few results have been obtained on some fractional integral inequalities using Saigo fractional integral operator, see \cite{C4,K4,P1,P2,YI}. Motivated from \cite{C1,C5,D1,KI2}, our purpose in this paper is to establish some new results using generalized k-fractional integral in terms of Gauss hypergeometric function. The paper has been organized as follows, in section 2, we define basic definitions and proposition related to generalized k-fractional integral. In section 3, we give the results about reverse Minkowski fractional integral inequality using fractional generalized k-fractional integral, In section 4, we give some other inequalities using fractional generalized k-fractional integral.
\section{Preliminaries}
\paragraph{} In this section, we give some necessary definitions which will be used latter.
\begin{definition} \cite{KI2,YI}
The function $f(x)$, for all $x>0$ is said to be in the $L_{p,k}[0,\infty),$ if
\begin{equation}
L_{p,k}[0,\infty)=\left\{f: \|f\|_{L_{p,k}[0,\infty)}=\left(\int_{0}^{\infty}|f(x)|^{p}x^{k}dx\right)^{\frac{1}{p}} < \infty \, \, 1 \leq p < \infty \, k \geq 0\right\},
\end{equation}
\end{definition}
\begin{definition} \cite{KI2,SAO,YI}
Let  $f \in L_{1,k}[0,\infty),$. The generalized Riemann-Liouville fractional integral $I^{\alpha,k}f(x)$ of order $\alpha, k \geq 0$ is defined by
\begin{equation}
I^{\alpha,k}f(x)= \frac{(k+1)^{1-\alpha}}{\Gamma (\alpha)}\int_{0}^{x}(x^{k+1}-t^{k+1})^{\alpha-1}t^{k} f(t)dt.
\end{equation}
\end{definition}
\begin{definition} \cite{KI2,YI}
Let $k\geq 0,\alpha>0, \mu >-1$ and $\beta, \eta \in R $. The generalized k-fractional integral $I^{\alpha,\beta,\eta,\mu}_{x,k}$ (in terms of the Gauss hypergeometric function)of order $\alpha$ for real-valued continuous function $f(t)$ is defined by
\begin{equation}
\begin{split}
I^{\alpha,\beta,\eta,\mu}_{x,k}[f(x)]&
= \frac{(k+1)^{\mu+\beta+1}x^{(k+1)(-\alpha-\beta-2\mu)}}{\Gamma (\alpha)}\int_{0}^{x}\tau^{(k+1)\mu}(x^{k+1}-\tau^{k+1})^{\alpha-1}
\times \\
& _{2}F_{1} (\alpha+ \beta+\mu, -\eta; \alpha; 1-(\frac{\tau}{x})^{k+1})\tau^{k} f(\tau)d\tau.
\end{split}
\end{equation}
\end{definition}
where, the function $_{2}F_{1}(-)$ in the right-hand side of (2.3) is the Gaussian hypergeometric function defined by
 \begin{equation}
 _{2}F_{1} (a, b; c; x)=\sum_{n=0}^{\infty}\frac{(a)_{n}(b)_{n}}{(c)_{n}} \frac{x^{n}}{n!},
\end{equation}
and $(a)_{n}$ is the Pochhammer symbol\\
$$(a)_{n}=a(a+1)...(a+n-1)=\frac{\Gamma(a+n)}{\Gamma(a)}, \,\,\,(a)_{0}=1.$$
Consider the function
\begin{equation}
\begin{split}
F(x,\tau)&= \frac{(k+1)^{\mu+\beta+1}x^{(k+1)(-\alpha-\beta-2\mu)}}{\Gamma (\alpha)}\tau^{(k+1)\mu}\\
&(x^{k+1}-\tau^{k+1})^{\alpha-1} \times _{2}F_{1} (\alpha+ \beta+\mu, -\eta; \alpha; 1-(\frac{\tau}{x})^{k+1})\\
&=\sum_{n=0}^{\infty}\frac{(\alpha+\beta+\mu)_{n}(-n)_{n}}{\Gamma(\alpha+n)n!}x^{(k+1)(-\alpha-\beta-2\mu-\eta)}\tau^{(k+1)\mu}(x^{k+1}-\tau^{k+1})^{\alpha-1+n}(k+1)^{\mu+\beta+1}\\
&=\frac{\tau^{(k+1)\mu}(x^{k+1}-\tau^{k+1})^{\alpha-1}(k+1)^{\mu+\beta+1}}{x^{k+1}(\alpha+\beta+2\mu)\Gamma(\alpha)}+\\
&\frac{\tau^{(k+1)\mu}(x^{k+1}-\tau^{k+1})^{\alpha}(k+1)^{\mu+\beta+1}(\alpha+\beta+\mu)(-n)}{x^{k+1}(\alpha+\beta+2\mu+1)\Gamma(\alpha+1)}+\\
&\frac{\tau^{(k+1)\mu}(x^{k+1}-\tau^{k+1})^{\alpha+1}(k+1)^{\mu+\beta+1}(\alpha+\beta+\mu)(\alpha+\beta+\mu+1)(-n)(-n+1)}{x^{k+1}(\alpha+\beta+2\mu+1)\Gamma(\alpha+2)2!}+...
\end{split}
\end{equation}
It is clear that $F(x,\tau)$ is positive  because for all $\tau \in (0, x)$ , $(x>0)$ since each term of the (2.5) is positive.
\section{Reverse Minkowski fractional integral inequality }
\paragraph{}In this section, we establish reverse Minkowski fractional integral inequality using generalized k-fractional integral operator (in terms of the Gauss hypergeometric function).
 \begin{theorem} Let $p\geq1$ and let $f$, $g$ be two positive function on $[0, \infty)$, such that for all $x>0$, $I^{\alpha,\beta,\eta,\mu}_{x,k}[f^{p}(x)]<\infty$, $I^{\alpha,\beta,\eta,\mu}_{x,k}[g^{q}(x)]<\infty$. If $0<m\leq \frac{f(\tau)}{g(\tau)}\leq M$, $\tau \in (0,x)$ we have
 \begin{equation}
\left[I^{\alpha,\beta,\eta,\mu}_{x,k}[f^{p}(x)]\right]^{\frac{1}{p}}+\left[I^{\alpha,\beta,\eta,\mu}_{x,k}[g^{q}(x)]\right]^{\frac{1}{p}}\leq \frac{1+M(m+2)}{(m+1)(M+1)}\left[I^{\alpha,\beta,\eta,\mu}_{x,k}[(f+g)^{p}(x)]\right]^{\frac{1}{p}},
\end{equation}
 for all $k \geq 0,$ $\alpha > max\{0,-\beta-\mu\}$, $\beta < 1,$ $\mu >-1,$ $\beta -1< \eta <0.$
 \end{theorem}
\textbf{Proof}: Using the condition $\frac{f(\tau)}{g(\tau)}\leq M$, $\tau \in (0,x)$, $x>0$, we can write
\begin{equation}
(M+1)^{p}f(\tau)\leq M^{p}(f+g)^{p}(\tau).
\end{equation}
Multiplying both side of (3.2) by $F(x,\tau)$, then integrating resulting identity with respect to $\tau$ from $0$ to $x$, we get
\begin{equation}
\begin{split}
&(M+1)^{p}\frac{(k+1)^{\mu+\beta+1}x^{(k+1)(-\alpha-\beta-2\mu)}}{\Gamma (\alpha)}\int_{0}^{x}\tau^{(k+1)\mu}(x^{k+1}-\tau^{k+1})^{\alpha-1}
\times \\
& _{2}F_{1} (\alpha+ \beta+\mu, -\eta; \alpha; 1-(\frac{\tau}{x})^{k+1})\tau^{k} f^{p}(\tau)d\tau\\
&\leq M^{p}\frac{(k+1)^{\mu+\beta+1}x^{(k+1)(-\alpha-\beta-2\mu)}}{\Gamma (\alpha)}\int_{0}^{x}\tau^{(k+1)\mu}(x^{k+1}-\tau^{k+1})^{\alpha-1}
\times \\
& _{2}F_{1} (\alpha+ \beta+\mu, -\eta; \alpha; 1-(\frac{\tau}{x})^{k+1})\tau^{k} (f+g)^{p}(\tau)d\tau,
\end{split}
\end{equation}
\noindent which is equivalent to
\begin{equation}
I^{\alpha,\beta,\eta,\mu}_{x,k}[f^{p}(x)] \leq \frac{M^{p}}{(M+1)^{p}} \left[I^{\alpha,\beta,\eta,\mu}_{x,k}[(f+g)^{p}(x)]\right],
\end{equation}
\noindent hence, we can write
\begin{equation}
\left[I^{\alpha,\beta,\eta,\mu}_{x,k}[f^{p}(x)] \right]^{\frac{1}{p}} \leq \frac{M}{(M+1)} \left[I^{\alpha,\beta,\eta,\mu}_{x,k}[(f+g)^{p}(x)]\right]^{\frac{1}{p}}.
\end{equation}
On other hand, using condition $m\leq \frac{f(\tau)}{g(\tau)}$, we obtain
\begin{equation}
(1+\frac{1}{m})g(\tau)\leq \frac{1}{m}(f(\tau)+g(\tau)),
\end{equation}
therefore,
\begin{equation}
(1+\frac{1}{m})^{p}g^{p}(\tau)\leq(\frac{1}{m})^{p}(f(\tau)+g(\tau))^{p}.
\end{equation}
Now, multiplying both side of (3.7) by  $F(x,\tau)$, ( $\tau \in(0,x)$, $x>0$), where  $G(x,\tau)$ is defined by (2.5). Then integrating resulting identity with respect to $\tau$ from $0$ to $x$, we have
\begin{equation}
\left[I^{\alpha,\beta,\eta,\mu}_{x,k}[g^{p}(x)]\right]^{\frac{1}{p}} \leq \frac{1}{(m+1)} \left[I^{\alpha,\beta,\eta,\mu}_{x,k}[(f+g)^{p}(x)]\right]^{\frac{1}{p}}.
\end{equation}
The inequalities (3.1) follows on adding the inequalities (3.5) and (3.8).
\paragraph{}Our second result is as follows.
\begin{theorem} Let $p\geq1$ and $f$, $g$ be two positive function on $[0, \infty)$, such that for all $x>0$, $I^{\alpha,\beta,\eta,\mu}_{x,k}[f^{p}(x)]<\infty$, $I^{\alpha,\beta,\eta,\mu}_{x,k}[g^{q}(x)]<\infty$. If $0<m\leq \frac{f(\tau)}{g(\tau)}\leq M$, $\tau \in (0,x)$  then we have
\begin{equation}
\begin{split}
\left[I^{\alpha,\beta,\eta,\mu}_{x,k}[f^{p}(x)] \right]^{\frac{2}{p}}+\left[I^{\alpha,\beta,\eta,\mu}_{x,k}[g^{q}(x)] \right]^{\frac{2}{p}}\geq &(\frac{(M+1)(m+1)}{M}-2)\left[I^{\alpha,\beta,\eta,\mu}_{x,k}[f^{p}(x)] \right]^{\frac{1}{p}}+\\
&\left[I^{\alpha,\beta,\eta,\mu}_{x,k}[g^{q}(x)] \right]^{\frac{1}{p}}.
\end{split}
\end{equation}
for all $k \geq 0,$ $\alpha > max\{0,-\beta-\mu\}$, $\beta < 1,$ $\mu >-1,$ $\beta -1< \eta <0.$
\end{theorem}
\textbf{Proof}: Multiplying the inequalities (3.5) and (3.8), we obtain
\begin{equation}
\frac{(M+1)(m+1)}{M}\left[I^{\alpha,\beta,\eta,\mu}_{x,k}[f^{p}(x)]\right]^{\frac{1}{p}}\times \left[I^{\alpha,\beta,\eta,\mu}_{x,k}[g^{q}(x)]\right]^{\frac{1}{p}}\leq \left([I^{\alpha,\beta,\eta,\mu}_{x,k}[(f(x)+g(x))^{p}]]^{\frac{1}{p}}\right)^{2}.
\end{equation}
Applying Minkowski inequalities to the right hand side of (3.10), we have
 \begin{equation}
 (\left[I^{\alpha,\beta,\eta,\mu}_{x,k}[(f(x)+g(x))^{p}]\right]^{\frac{1}{p}})^{2}\leq (\left[I^{\alpha,\beta,\eta,\mu}_{x,k}[f^{p}(x)]\right]^{\frac{1}{p}}+\left[I^{\alpha,\beta,\eta,\mu}_{x,k}[g^{q}(x)]\right]^{\frac{1}{p}})^{2},
\end{equation}
which implies that
\begin{equation}
\begin{split}
 \left[I^{\alpha,\beta,\eta,\mu}_{x,k}[(f(x)+g(x))^{p}]\right]^{\frac{2}{p}}\leq & \left[I^{\alpha,\beta,\eta,\mu}_{x,k}[f^{p}(x)]\right]^{\frac{2}{p}}+
 \left[I^{\alpha,\beta,\eta,\mu}_{x,k}[g^{q}(x)]\right]^{\frac{2}{p}}\\
 &+2\left[I^{\alpha,\beta,\eta,\mu}_{x,k}[f^{p}(x)]\right]^{\frac{1}{p}} \left[I^{\alpha,\beta,\eta,\mu}_{x,k}[g^{q}(x)]\right]^{\frac{1}{p}}.
\end{split}
\end{equation}
Hence, from (3.10) and (3.12), we obtain (3.9).
Theorem 3.2 is thus proved.
\section{ Other fractional integral inequalities related to Minkowski inequality}
\paragraph{}In this section, we establish some new integral inequalities related to Minkowski inequality using generalized k-fractional integral operator (in terms of the Gauss hypergeometric function).
\begin{theorem} Let $p>1$,  $\frac{1}{p}+\frac{1}{q}=1 $ and $f$, $g$ be two positive function on $[0, \infty)$, such that $I^{\alpha,\beta,\eta,\mu}_{x,k}[f(x)]<\infty$, $I^{\alpha,\beta,\eta,\mu}_{x,k}[g(x)]<\infty$. If $0<m\leq \frac{f(\tau)}{g(\tau)}\leq M < \infty$, $\tau \in [0,x]$ we have
\begin{equation}
\left[I^{\alpha,\beta,\eta,\mu}_{x,k}[f(x)]\right]^{\frac{1}{p}} \left[I^{\alpha,\beta,\eta,\mu}_{x,k}[g(x)]\right]^{\frac{1}{q}}
\leq (\frac{M}{m})^{\frac{1}{pq}}\left[I^{\alpha,\beta,\eta,\mu}_{x,k}[[f(x)]^{\frac{1}{p}}[g(x)]^{\frac{1}{q}}]\right],
\end{equation}
 for all $k \geq 0,$ $\alpha > max\{0,-\beta-\mu\}$, $\beta < 1,$ $\mu >-1,$ $\beta -1< \eta <0.$
\end{theorem}
\textbf{Proof:-} Since $\frac{f(\tau)}{g(\tau)}\leq M $, $\tau \in[0,x]$  $x> 0$, therefore \\
\begin{equation}
[g(\tau)]^{\frac{1}{p}}\geq M^{\frac{-1}{q}}[f(\tau)]^{\frac{1}{q}},
\end{equation}
and also,
\begin{equation}
\begin{split}
[f(\tau)]^{\frac{1}{p}}[g(\tau)]^{\frac{1}{q}}&\geq M^{\frac{-1}{q}}[f(\tau)]^{\frac{1}{q}}[f(\tau)]^{\frac{1}{p}}\\
&\geq  M^{\frac{-1}{q}}[f(\tau)]^{\frac{1}{q}+\frac{1}{q}}\\
&\geq  M^{\frac{-1}{q}}[f(\tau)].
\end{split}
\end{equation}
Multiplying both side of (4.3) by $F(x,\tau)$, ( $\tau \in(0,x)$, $x>0$), where  $F(x,\tau)$ is defined by (2.5). Then integrating resulting identity with respect to $\tau$ from $0$ to $x$, we have
\begin{equation}
\begin{split}
&\frac{(k+1)^{\mu+\beta+1}x^{(k+1)(-\alpha-\beta-2\mu)}}{\Gamma (\alpha)}\int_{0}^{x}\tau^{(k+1)\mu}(x^{k+1}-\tau^{k+1})^{\alpha-1}
\times \\
& _{2}F_{1} (\alpha+ \beta+\mu, -\eta; \alpha; 1-(\frac{\tau}{x})^{k+1})\tau^{k} f(\tau)^{\frac{1}{p}}g(\tau)^{\frac{1}{q}}d\tau \\
&\leq M^{\frac{-1}{q}}\frac{(k+1)^{\mu+\beta+1}x^{(k+1)(-\alpha-\beta-2\mu)}}{\Gamma (\alpha)}\int_{0}^{x}\tau^{(k+1)\mu}(x^{k+1}-\tau^{k+1})^{\alpha-1}
\times \\
& _{2}F_{1} (\alpha+ \beta+\mu, -\eta; \alpha; 1-(\frac{\tau}{x})^{k+1})\tau^{k} f(\tau)d\tau,
\end{split}
\end{equation}
which implies that
\begin{equation}
I^{\alpha,\beta,\eta,\mu}_{x,k}\left[[f(x)]^{\frac{1}{p}}[g(x)]^{\frac{1}{q}} \right] \leq M^{\frac{-1}{q}} \left[I^{\alpha,\beta,\eta,\mu}_{x,k}f(x)\right].
\end{equation}
Consequently,
\begin{equation}
\left(I^{\alpha,\beta,\eta,\mu}_{x,k}\left[[f(x)]^{\frac{1}{p}}[g(x)]^{\frac{1}{q}} \right]\right)^{\frac{1}{p}} \leq M^{\frac{-1}{pq}} \left[I^{\alpha,\beta,\eta,\mu}_{x,k}f(x)\right]^{\frac{1}{p}},
\end{equation}
on other hand, since $m g(\tau)\leq f(\tau)$, \, $\tau \in[0,x)$, $x>0$, then we have
\begin{equation}
[f(\tau)]^{\frac{1}{p}}\geq m^{\frac{1}{p}}[g(\tau)]^{\frac{1}{p}},
\end{equation}
multiplying equation (4.7) by $[g(\tau)]^{\frac{1}{q}}$, we have
\begin{equation}
[f(\tau)]^{\frac{1}{p}}[g(\tau)]^{\frac{1}{q}}\geq m^{\frac{1}{p}}[g(\tau)]^{\frac{1}{q}}[g(\tau)]^{\frac{1}{p}}= m^{\frac{1}{p}}[g(\tau)].
\end{equation}
Multiplying both side of (4.8) by $F(x,\tau)$, ( $\tau \in(0,x)$, $x>0$), where  $F(x,\tau)$ is defined by (2.5). Then integrating resulting identity with respect to $\tau$ from $0$ to $x$, we have
\begin{equation}
\begin{split}
&\frac{(k+1)^{\mu+\beta+1}x^{(k+1)(-\alpha-\beta-2\mu)}}{\Gamma (\alpha)}\int_{0}^{x}\tau^{(k+1)\mu}(x^{k+1}-\tau^{k+1})^{\alpha-1}
\times \\
& _{2}F_{1} (\alpha+ \beta+\mu, -\eta; \alpha; 1-(\frac{\tau}{x})^{k+1})\tau^{k} f(\tau)^{\frac{1}{p}}g(\tau)^{\frac{1}{q}}d\tau  \\
&\leq  M^{\frac{1}{p}}\frac{(k+1)^{\mu+\beta+1}x^{(k+1)(-\alpha-\beta-2\mu)}}{\Gamma (\alpha)}\int_{0}^{x}\tau^{(k+1)\mu}(x^{k+1}-\tau^{k+1})^{\alpha-1}
\times \\
& _{2}F_{1} (\alpha+ \beta+\mu, -\eta; \alpha; 1-(\frac{\tau}{x})^{k+1})\tau^{k} g(\tau)d\tau,
\end{split}
\end{equation}
 which implies that
\begin{equation}
I^{\alpha,\beta,\eta,\mu}_{x,k}\left[[f(x)]^{\frac{1}{p}}[g(x)]^{\frac{1}{q}} \right] \leq m^{\frac{1}{p}} \left[I^{\alpha,\beta,\eta,\mu}_{x,k}g(x)\right].
\end{equation}
Hence, we can write
\begin{equation}
\left(I^{\alpha,\beta,\eta,\mu}_{x,k}\left[[f(x)]^{\frac{1}{p}}[g(x)]^{\frac{1}{q}} \right]\right)^{\frac{1}{q}} \leq m^{\frac{1}{pq}} \left[I^{\alpha,\beta,\eta,\mu}_{x,k}f(x)\right]^{\frac{1}{q}},
\end{equation}
multiplying equation (4.6) and (4.11) we get the result (4.1).
\begin{theorem} Let $f$ and $g$ be two positive function on $[0, \infty[$, such that\\  $I^{\alpha,\beta,\eta,\mu}_{x,k}[f^{p}(x)]<\infty$,
 $I^{\alpha,\beta,\eta,\mu}_{x,k}[g^{q}(x)]<\infty$. $x>0$,  If $0<m\leq \frac{f(\tau)^{p}}{g(\tau)^{q}}\leq M < \infty$, $\tau \in [0,x]$. Then we have
 \begin{equation*}
\left[I^{\alpha,\beta,\eta,\mu}_{x,k}f^{p}(x)\right]^{\frac{1}{p}} \left[I^{\alpha,\beta,\eta,\mu}_{x,k}g^{q}(x)\right]^{\frac{1}{q}}\leq (\frac{M}{m})^{\frac{1}{pq}}\left[I^{\alpha,\beta,\eta,\mu}_{x,k}(f(x)g(x))\right] hold.
\end{equation*}
Where  $p>1$,  $\frac{1}{p}+\frac{1}{q}=1 $, for all $k \geq 0,$ $\alpha > max\{0,-\beta-\mu\}$, $\beta < 1,$ $\mu >-1,$ $\beta -1< \eta <0.$
\end{theorem}
\textbf{Proof:-}
Replacing $f(\tau)$ and $g(\tau)$ by $f(\tau)^{p}$ and $g(\tau)^{q}$,  $\tau \in [0,x]$, $x>0$ in theorem 4.1, we obtain required inequality.
\paragraph{} Now, here we present fractional integral inequality related to Minkowsky inequality as follows
\begin{theorem} let $f$ and $g$ be two integrable functions on $[1, \infty]$ such that $\frac{1}{p}+\frac{1}{q}=1, p>1,$ and $0<m<\frac{f(\tau)}{g(\tau)}<M, \tau \in (0,x), x>0.$ Then for all $\alpha>0,$ we have
\begin{equation}
I^{\alpha,\beta,\eta,\mu}_{x,k}\{fg\}(x)\leq \frac{2^{p-1}M^{p}}{p(M+1)^{p}}\left(I^{\alpha,\beta,\eta,\mu}_{x,k}[f^{p}+g^{p}](x)\right)+\frac{2^{q-1}}{q(m+1)^{q}}\left(I^{\alpha,\beta,\eta,\mu}_{x,k}[f^{q}+g^{q}](x)\right),
\end{equation}
 for all $k \geq 0,$ $\alpha > max\{0,-\beta-\mu\}$, $\beta < 1,$ $\mu >-1,$ $\beta -1< \eta <0.$
\end{theorem}
\textbf{Proof:-} Since,  $\frac{f(\tau)}{g(\tau)}<M, \tau \in (0,x), x>0,$ we have
\begin{equation}
(M+1)f(\tau)\leq M(f+g)(\tau).
\end{equation}
Taking $p^{th}$  power on both side and multiplying resulting identity by $ F(x,\tau)$, we obtain
\begin{equation}
\begin{split}
&(M+1)^{p}\frac{(k+1)^{\mu+\beta+1}x^{(k+1)(-\alpha-\beta-2\mu)}}{\Gamma (\alpha)}\int_{0}^{x}\tau^{(k+1)\mu}(x^{k+1}-\tau^{k+1})^{\alpha-1}
\times \\
& _{2}F_{1} (\alpha+ \beta+\mu, -\eta; \alpha; 1-(\frac{\tau}{x})^{k+1})\tau^{k} f^{p}(\tau)d\tau\\
 &\leq M^{p} \frac{(k+1)^{\mu+\beta+1}x^{(k+1)(-\alpha-\beta-2\mu)}}{\Gamma (\alpha)}\int_{0}^{x}\tau^{(k+1)\mu}(x^{k+1}-\tau^{k+1})^{\alpha-1}
\times \\
& _{2}F_{1} (\alpha+ \beta+\mu, -\eta; \alpha; 1-(\frac{\tau}{x})^{k+1})\tau^{k} (f+g)^{p}(\tau)d\tau,
\end{split}
\end{equation}
therefore,
\begin{equation}
I^{\alpha,\beta,\eta,\mu}_{x,k}[f^{p}(x)]\leq \frac{M^{p}}{(M+1)^{p}}I^{\alpha,\beta,\eta,\mu}_{x,k}[(f+g)^{p}(x)],
\end{equation}
on other hand, $0<m<\frac{f(\tau)}{g(\tau)}, \tau \in (0,x), x>0,$ we can write
\begin{equation}
(m+1)g(\tau)\leq (f+g)(\tau),
\end{equation}
therefore,
\begin{equation}
\begin{split}
&(m+1)^{q}\frac{(k+1)^{\mu+\beta+1}x^{(k+1)(-\alpha-\beta-2\mu)}}{\Gamma (\alpha)}\int_{0}^{x}\tau^{(k+1)\mu}(x^{k+1}-\tau^{k+1})^{\alpha-1}
\times \\
& _{2}F_{1} (\alpha+ \beta+\mu, -\eta; \alpha; 1-(\frac{\tau}{x})^{k+1})\tau^{k} g^{q}(\tau)d\tau\\
 &\leq  \frac{(k+1)^{\mu+\beta+1}x^{(k+1)(-\alpha-\beta-2\mu)}}{\Gamma (\alpha)}\int_{0}^{x}\tau^{(k+1)\mu}(x^{k+1}-\tau^{k+1})^{\alpha-1}
\times \\
& _{2}F_{1} (\alpha+ \beta+\mu, -\eta; \alpha; 1-(\frac{\tau}{x})^{k+1})\tau^{k} (f+g)^{q}(\tau)d\tau,
\end{split}
\end{equation}
consequently, we have
\begin{equation}
I^{\alpha,\beta,\eta,\mu}_{x,k}[g^{q}(x)]\leq \frac{1}{(m+1)^{q}}I^{\alpha,\beta,\eta,\mu}_{x,k}[(f+g)^{q}(x)].
\end{equation}
Now, using Young inequality
\begin{equation}
[f(\tau)g(\tau)]\leq \frac{f^{p}(\tau)}{p}+\frac{g^{q}(\tau)}{q}.
\end{equation}
Multiplying both side of (4.19) by $ F(x,\tau)$, which is positive because $\tau \in(0,x)$, $x>0$, then integrate the resulting identity with respect to $\tau$ from $0$ to $x$, we get
\begin{equation}
I^{\alpha,\beta,\eta,\mu}_{x,k}[f(x)g(x))]\leq \frac{1}{p}\,I^{\alpha,\beta,\eta,\mu}_{x,k}[f^{p}(x)]+\frac{1}{q}\,I^{\alpha,\beta,\eta,\mu}_{x,k}[g^{q}(x)],
\end{equation}
from equation (4.15), (4.18) and (4.20) we get
\begin{equation}
I^{\alpha,\beta,\eta,\mu}_{x,k}[f(x)g(x))]\leq \frac{M^{p}}{p(M+1)^{p}}\,I^{\alpha,\beta,\eta,\mu}_{x,k}[(f+g)^{p}(x)]+\frac{1}{q(m+1)^{q}}\,I^{\alpha,\beta,\eta,\mu}_{x,k}[(f+g)^{q}(x)],
\end{equation}
now using the inequality $(a+b)^{r}\leq 2^{r-1}(a^{r}+b^{r}), r>1, a,b \geq 0,$ we have
\begin{equation}
I^{\alpha,\beta,\eta,\mu}_{x,k}[(f+g)^{p}(x)] \leq 2^{p-1}I^{\alpha,\beta,\eta,\mu}_{x,k}[(f^{p}+g^{p})(x)],
\end{equation}
and
\begin{equation}
I^{\alpha,\beta,\eta,\mu}_{x,k}[(f+g)^{q}(x)] \leq 2^{q-1}I^{\alpha,\beta,\eta,\mu}_{x,k}[(f^{q}+g^{q})(x)].
\end{equation}
Injecting (4.22), (2.23) in (4.21) we get required inequality (4.12). This complete the proof.
\begin{theorem}
Let $f$, $g$ be two positive function on $[0, \infty)$, such that $f$ is non-decreasing and $g$ is non-increasing. Then
\begin{equation}
\begin{split}
I^{\alpha,\beta,\eta,\mu}_{x,k}f^{\gamma}(x) g^{\delta}(x)&\leq (k+1)^{-\mu-\beta}x^{(k+1)(\mu+\beta)}\frac{\Gamma(1-\beta)\Gamma(1+\mu+\eta+1)}{\Gamma(1-\beta+\eta)\Gamma(\mu+1)} \\
&\times  I^{\alpha,\beta,\eta,\mu}_{x,k}[f^{\gamma}(x)]I^{\alpha,\beta,\eta,\mu}_{x,k}[g^{\delta}(x)],
\end{split}\end{equation}
 for all $k \geq 0,$ $\alpha > max\{0,-\beta-\mu\}$, $\beta < 1,$ $\mu >-1,$ $\beta -1< \eta <0.$
\end{theorem}
\textbf{Proof:-} let $\tau,\rho \in [0,x]$, $x>0$, for any $\delta>0$, $\gamma>0$, we have
\begin{equation}
\left(f^{\gamma}(\tau)-f^{\gamma}(\rho)\right)\left(g^{\delta}(\rho)-g^{\delta}(\tau)\right) \geq 0,
\end{equation}
\begin{equation}
f^{\gamma}(\tau)g^{\delta}(\rho)-f^{\gamma}(\tau)g^{\delta}(\tau)- f^{\gamma}(\rho)(g^{\delta}(\rho)+f^{\gamma}(\rho)g^{\delta}(\tau) \geq 0,
\end{equation}
therefore
\begin{equation}
f^{\gamma}(\tau)g^{\delta}(\tau)+f^{\gamma}(\rho)(g^{\delta}(\rho)\leq f^{\gamma}(\tau)g^{\delta}(\rho)+f^{\gamma}(\rho)g^{\delta}(\tau).
\end{equation}
Now, multiplying both side of (4.27) by $F(x,\tau)$, ( $\tau \in(0,x)$, $x>0$), where  $F(x,\tau)$ is defined by (2.5). Then integrating resulting identity with respect to $\tau$ from $0$ to $x$, we have
\begin{equation}
\begin{split}
&\frac{(k+1)^{\mu+\beta+1}x^{(k+1)(-\alpha-\beta-2\mu)}}{\Gamma (\alpha)}\int_{0}^{x}\tau^{(k+1)\mu}(x^{k+1}-\tau^{k+1})^{\alpha-1}
\times \\
& _{2}F_{1} (\alpha+ \beta+\mu, -\eta; \alpha; 1-(\frac{\tau}{x})^{k+1})\tau^{k}[f^{\gamma}(\tau)g^{\delta}(\tau)]d\tau\\
&+ f^{\gamma}(\rho)g^{\delta}(\rho)\frac{(k+1)^{\mu+\beta+1}x^{(k+1)(-\alpha-\beta-2\mu)}}{\Gamma (\alpha)}\int_{0}^{x}\tau^{(k+1)\mu}(x^{k+1}-\tau^{k+1})^{\alpha-1}
\times \\
& _{2}F_{1} (\alpha+ \beta+\mu, -\eta; \alpha; 1-(\frac{\tau}{x})^{k+1})\tau^{k}[1]d\tau \\
&\leq g^{\delta}(\rho)\frac{(k+1)^{\mu+\beta+1}x^{(k+1)(-\alpha-\beta-2\mu)}}{\Gamma (\alpha)}\int_{0}^{x}\tau^{(k+1)\mu}(x^{k+1}-\tau^{k+1})^{\alpha-1}
\times \\
& _{2}F_{1} (\alpha+ \beta+\mu, -\eta; \alpha; 1-(\frac{\tau}{x})^{k+1})\tau^{k}f^{\gamma}(\tau)d\tau\\
&+f^{\gamma}(x)\frac{(k+1)^{\mu+\beta+1}x^{(k+1)(-\alpha-\beta-2\mu)}}{\Gamma (\alpha)}\int_{0}^{x}\tau^{(k+1)\mu}(x^{k+1}-\tau^{k+1})^{\alpha-1}
\times \\
& _{2}F_{1} (\alpha+ \beta+\mu, -\eta; \alpha; 1-(\frac{\tau}{x})^{k+1})\tau^{k} g^{\delta}(\tau)d\tau,
\end{split}
\end{equation}
\begin{equation}
\begin{split}
&I^{\alpha,\beta,\eta,\mu}_{x,k}[f^{\gamma}(x)g^{\delta}(x)]+f^{\gamma}(\rho)(g^{\delta}(\rho)I^{\alpha,\beta,\eta,\mu}_{x,k}[1]\\
&\leq g^{\delta}(\rho)I^{\alpha,\beta,\eta,\mu}_{x,k}[f^{\gamma}(x)]+f^{\gamma}(\rho)I^{\alpha,\beta,\eta,\mu}_{x,k}[g^{\delta}(x)].
\end{split}
\end{equation}
Again, multiplying both side of (4.29) by  $F(x,\rho)$, ( $\rho \in(0,x)$, $x>0$), where  $F(x,\rho)$ is defined by (2.5). Then integrating resulting identity with respect to $\rho$ from $0$ to $x$, we have
\begin{equation*}
\begin{split}
&I^{\alpha,\beta,\eta,\mu}_{x,k}[f^{\gamma}(x)g^{\delta}(x)]I^{\alpha,\beta,\eta,\mu}_{x,k}[1]+I^{\alpha,\beta,\eta,\mu}_{x,k}[f^{\gamma}(x)g^{\delta}(x)]
I^{\alpha,\beta,\eta,\mu}_{x,k}[1]\\
&\leq I^{\alpha,\beta,\eta,\mu}_{x,k}[g^{\delta}(x)]I^{\alpha,\beta,\eta,\mu}_{x,k}[f^{\gamma}(x)]
+I^{\alpha,\beta,\eta,\mu}_{x,k}[f^{\gamma}(x)]I^{\alpha,\beta,\eta,\mu}_{x,k}[g^{\delta}(x)],
\end{split}
\end{equation*}
then we can write
\begin{equation*}
2I^{\alpha,\beta,\eta,\mu}_{x,k}[f^{\gamma}(x)g^{\delta}(x)] \leq \frac{1}{[I^{\alpha,\beta,\eta,\mu}_{x,k}[1]]^{-1}}2 II^{\alpha,\beta,\eta,\mu}_{x,k}[f^{\gamma}(x)]I^{\alpha,\beta,\eta,\mu}_{x,k}[g^{\delta}(x)].
\end{equation*}
This proves the result (4.24).\\
\textbf{Competing interests}\\
The authors declare that they have no competing interests.
		
\end{document}